\documentclass[11pt]{amsart}
\usepackage{epsfig}
\usepackage{amsmath}
\usepackage{amsthm}
\newtheorem{theorem}{Theorem}[section]
\newtheorem{lemma}[theorem]{Lemma}
\newtheorem{proposition}[theorem]{Proposition}

\theoremstyle{definition}

\begin{document}

 \title
{A family of bijections between $G$-parking functions and spanning trees}

  \author {Denis Chebikin and Pavlo Pylyavskyy}
\address{Department of Mathematics, Massachusetts Institute of
  Technology,
Cambridge, MA 02139}
\email{chebikin@mit.edu, pasha@mit.edu}

\begin{abstract}
For a directed graph $G$ on vertices $\{0,1,\dots,n\}$, a
\emph{$G$-parking function} is an $n$-tuple $(b_1,\dots,b_n)$
of non-negative integers such that, for every non-empty subset
$U\subseteq \{1,\dots,n\}$, there exists a vertex $j\in U$ for which
there are more than $b_j$ edges going from $j$ to $G-U$. We construct
a family of bijective maps between the set $\mathcal P_G$ of $G$-parking
functions and the set $\mathcal T_G$ of spanning trees of $G$ rooted
at $0$, thus providing a combinatorial proof of $|\mathcal P_G| =
|\mathcal T_G|$. 
\end{abstract}

  \maketitle 

\section {Introduction}

%Let $G$ be a directed graph on $n+1$ vertices indexed by integers from
%$0$ to $n$. We allow $G$ to have multiple edges but not loops. A
%{\it spanning tree} of $G$ rooted at $0$ is a subgraph of $G$ such
%that, 
%for each $i
%\in \{1, 2, \dotsm, n\}$, there is a unique path from $i$ to
%$0$ along the edges of the spanning tree. Note that these are the usual
%spanning trees of the graph with each edge
%oriented towards $0$. The number of such trees is given
%by the Matrix-Tree Theorem, see [Stan]. For the complete graph
%$K_{n+1}$, this number is $(n+1)^{n-1}$ by Cayley theorem. 
%
%We also consider $G$-parking functions, introduced in [Post]. A
%\emph{$G$-parking function} is a sequence $(b_1, \dotsm,
%b_n)$ of non-negative integers 
%that satisfies the following condition: for each subset $U
%\subseteq \{1, 2, \dotsm, n\}$
%of
%vertices of $G$, there exists a vertex $j
%\in U$ such that the number of edges from $j$ to vertices outside of 
%$U$ is greater than $b_j$. In [Post], it
%is shown that the number of spanning trees of $G$
%is equal to the number of $G$-parking functions for any digraph $G$. In this
%paper, we construct two bijective maps bewteen the set $\mathcal P_G$
%of $G$-parking functions and the set $\mathcal T_G$ of spanning trees of $G$.
%
%Throughout the paper all spanning trees are rooted at 0 unless
%otherwise specified.

The classical parking functions are defined in the following way. There are $n$
drivers,
labeled $1,\dots,n$,  
and $n$ parking spots, $0,\dots, n-1$, arranged linearly in this
order. 
Each driver $i$ has a favorite parking spot $b_i$. Drivers enter 
the parking area in the order in which they are labeled. Each
driver proceeds to his favorite spot and parks 
there if it is free, or parks at the next available spot otherwise. 
The sequence $(b_1,\dots,b_n)$ is called 
a {\it parking function} if every driver parks successfully by this
rule. The most notable result about parking 
functions is a
bijective correspondence between such functions and trees on $n+1$ labeled
vertices. The number of such trees is $(n+1)^{n-1}$ by Cayley's
theorem. 
For more on 
parking functions, see for example \cite{Stan}.

Postnikov and Shapiro \cite{Post} suggested the following generalization
of parking functions. Let $G$ be a directed graph on $n+1$ vertices indexed 
by integers from
$0$ to $n$. A
\emph{$G$-parking function} is a sequence $(b_1, \dots,
b_n)$ of non-negative integers
that satisfies the following condition: for each subset $U
\subseteq \{1, 2, \dots, n\}$
of
vertices of $G$, there exists a vertex $j
\in U$ such that the number of edges from $j$ to vertices outside of
$U$ is greater than $b_j$. For the complete graph $G = K_{n+1}$, these
are the classical parking functions (we view $K_{n+1}$ as the digraph
with exactly one edge $(i,j)$ for all $i\neq j$).

A {\it spanning tree} of $G$ rooted at $m$ is a subgraph of $G$ such
that,
for each $i
\in \{0, 1, \dotsm, n\}$, there is a unique path from $i$ to
$m$ along the edges of the spanning tree. Note that these are the
spanning trees of the graph in the usual sense with each edge
oriented towards $m$. The number of such trees is given
by the Matrix-Tree Theorem; see \cite{Stan}. In \cite{Post} it
is shown that the number of spanning trees of $G$ rooted at $0$
is equal to the number of $G$-parking functions for any digraph $G$.
 
An equivalent fact was originally discovered by Dhar \cite{Dhar}, 
who studied the sandpile 
model. The so called recurrent states of the sandpile model are in 
one-to-one correspondence with $G$-parking functions for certain graphs 
$G$, including all symmetric graphs. 
A bijection between recurrent states and spanning trees for symmetric 
graphs $G$ is mentioned in \cite{Iva}, and a class of bijections is
constructed in \cite{Cori}.
The sandpile model was also studied by Gabrielov in \cite{Gab}. This paper
also
contains an extensive list of references on the topic.

In this paper
we present a family of bijections between $G$-parking functions and rooted 
spanning trees
of $G$. Given a spanning tree $T$ of $G$, we establish a total order
on the vertices of $T$ satisfying two conditions, and each such order
gives rise to a bijection in the family. In \cite{Fran} Francon used a similar
concept, which he called selection procedures, to construct a family
of bijections between parking functions and rooted trees in the
classical case $G=K_{n+1}$. Thus our result provides a generalization
of Francon's construction.

\section {A family of bijections}\label{bijections}

Let $G$ be a directed graph on vertices $\{0,\dots,n\}$. We allow $G$
to have multiple edges but not loops. To distinguish between multiple
edges of $G$, we fix an order on the set of edges going from $i$ to
$j$ for all $i\neq j$.

A \emph{subtree} of $G$ rooted at $m$ is a subgraph $T$ of $G$ containing
$m$ such that for every vertex $i$ of $T$, there is a unique path in
$T$ from $i$ to $m$. A subtree is called a \emph{spanning tree} if it
contains all vertices of $G$. 

Let $\mathbb T_G$ be the set of subtrees of $G$ rooted at $0$, and
let $\mathcal T_G$ be the set of spanning trees of $G$ rooted at
$0$. Unless stated otherwise, all spanning trees in this paper are
assumed to be rooted at $0$. Let $\mathcal P_G$ be the set of
$G$-parking functions. In this section we give a
bijection between $\mathcal T_G$ and $\mathcal P_G$.

For every $T\in \mathbb T_G$, let $\pi(T)$ be a total
order on the vertices of $T$, and write $i <_{\pi(T)} j$ to denote that
$i$ is smaller than $j$ in this order. We call the set $\Pi(G)=\{\pi(T)\ |\
T\in\mathbb T_G\}$ a \emph{proper set of tree orders} if the following
conditions hold for all $T\in \mathbb{T}_G$:
\begin{enumerate}

\item if $(j,i)$ is an edge of $T$, then $i<_{\pi(T)} j$;

\item if $t$ is a subtree of $T$ rooted at $0$, then the order
  $\pi(t)$ is consistent with $\pi(T)$; in other words, $i <_{\pi(t)}
  j$ if and only if $i <_{\pi(T)} j$ for $i,j\in t$.

\end{enumerate}
We give several examples of proper sets of tree orders in Section
\ref{examples}.

For $T\in \mathbb T_G$ and a vertex $j$ of $G$, the order $\pi(T)$
induces the
order on the edges going from $j$ to vertices of $T$ in which $(j,i)$
is smaller than $(j,i')$ whenever $i <_{\pi(T)} i'$ and which is
consistent with the previously fixed order on multiple edges. We write
$e <_{\pi(T)} e'$ to denote that $e$ is smaller than $e'$ in this order.

Given a proper set of tree orders $\Pi(G)$, define the map $\Theta_{\Pi,G} :
\mathcal T_G \rightarrow \mathcal P_G$ as follows. For $T\in \mathcal
T_G$ and a vertex $j\in \{1,\dots,n\}$, let $e_j$ be the edge of $T$
going out of $j$. Set $\Theta_{\Pi,G}(T) = (b_1,\dots,b_n)$, where $b_j$ is the
number of edges $e$ going out of $j$ such that $e<_{\pi(T)} e_j$. For
the rest of the section, we write $\Theta$ instead of $\Theta_{\Pi,G}$.

\begin{theorem}\label{main}
The map $\Theta$ is a bijection between $\mathcal T_G$ and $\mathcal P_G$.
\end{theorem}

\begin{proof}
We begin by checking that $\Theta(T)$ is a $G$-parking function.

\begin{lemma}\label{isPF}
$\Theta(T) \in \mathcal P_G$ for $T\in \mathcal T_G$.
\end{lemma}

\begin{proof}
For a subset $U\subseteq \{1,\dots,n\}$, let $j$ be the smallest
vertex of $U$ in the order $\pi(T)$. Let $e_j=(j,i)$ be the edge of
$T$ coming out of $j$. Then $i <_{\pi(T)} j$, so $i\notin U$ by choice
of $j$. For each of the $b_j$ edges 
$e=(j,i')$ such that $e<_{\pi(T)} e_j$, we have
$i' \leq_{\pi(T)} i <_{\pi(T)} j$, so $i'\notin U$. Thus there are at least
$b_j+1$ edges going from $j$ to vertices outside of $U$.
\end{proof}

Next, we define the inverse map $\Phi_{\Pi,G} : \mathcal P_G \rightarrow
\mathcal
T_G$. Given $P=(b_1,\dots,b_n)
\in \mathcal P_G$, we construct the corresponding tree 
$\Phi_{\Pi,G}(P)$ one edge at a time. Initially, let $t_0$ be the subtree of
$G$ consisting of the vertex $0$ alone, and put $p_0=0$. For $1\leq m
\leq n$, we choose the vertex $p_m$ and construct the subtree $t_m$
rooted at $0$ inductively as follows. Let $U_m$ be the set of vertices
not in $t_{m-1}$, and let $V_m$ be the set of vertices
$j\in U_m$ such that the number of edges from $j$ to $t_{m-1}$ is at
least $b_j+1$. Note that $|V_m|\geq 1$ by definition of a $G$-parking function.
For each $j\in V_m$, let $e_j$ be the edge from $j$ to $t_{m-1}$ such
that exactly $b_j$ edges $e$ from $j$ to $t_{m-1}$ satisfy
$e<_{\pi(t_{m-1})} e_j$. Let $t$ be the tree obtained by adjoining
each vertex $j\in V_m$ to $t_{m-1}$ by means of the edge $e_j$. Set
$p_m$ to be the smallest vertex of $V_m$ in the order $\pi(t)$, and set
$t_m$ to be the tree obtained by adjoining $p_m$ to $t_{m-1}$ by means
of the edge $e_{p_m}$. 
Obviously, $t_m$ is a subtree of $G$.
In the end, set $\Phi_{\Pi,G}(P) = T=t_n$. For
the rest of the section, we write $\Phi$ instead of $\Phi_{\Pi,G}$.

\begin{figure}
\begin{center}
\input{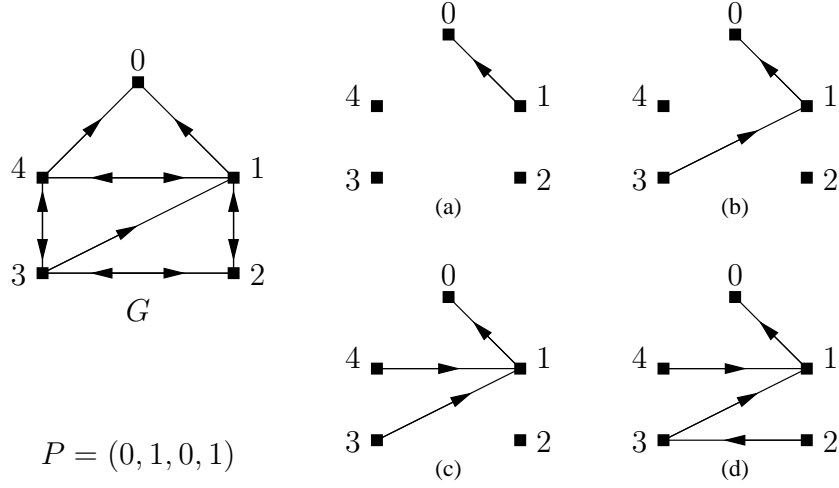}
\end{center}
\caption{An example of constructing $\Phi(P)$.}\label{exempel2}
\end{figure}

An example of constructing $\Phi(P)$ is shown in Figure
\ref{exempel2}. 
Let $G$ be the graph shown in the figure, and let $P=(0,1,0,1)$.
Let $\Pi$ be the tree order in which vertex $i$ is
smaller than vertex $j$ if $i$ is closer to the root than $j$, or else
if $i$ and $j$ are equidistant to the root, and $i<j$. Initially, $U_1 = \{1,2,3,4\}$ and $V_1 = \{1\}$, so
vertex $1$ is attached to the root to produce the subtree $t_1$. Then
we have $V_2 = \{3,4\}$ with $e_3 = (3,1)$ and $e_4 =
(4,1)$. Adjoining vertices $3$ and $4$ to $t_1$ by means of $e_3$ and
$e_4$ places vertices $3$ and $4$ the same distance away from the
root, making $p_2 = 3$,
so vertex $3$ is attached by means of the edge
$e_3 = (3,1)$ to produce $t_2$. At the next step, we have $V_3 = \{2,4\}$ with $e_2 =
(2,3)$ and $e_4 = (4,1)$. Adjoining vertices $2$ and $4$ to $t_2$ by
means of $e_2$ and $e_4$ makes vertex $4$ closer to the root than
vertex $2$, so we select vertex $4$ and attach it to vertex $1$ to
form $t_3$. Finally, we attach vertex $2$ to vertex $3$ to form $t_4 = \Phi(P)$.

%\begin{lemma}\label{numedges}
%For a vertex $k\in \{1,\dots,n\}$ and the edge $e_k$ coming out of $k$
%in $T$, the number of edges $e$ coming out of $k$ such that
%$e<_{\pi(T)} e_k$ is at least $b_k$.
%\end{lemma}

%\begin{proof}
%We have $k=p_m$ for some $1\leq m \leq n$. Since $k\in V$ at the step
%when $t_m$ is constructed, there are exactly $b_k$ edges $e$ going from
%$k$ to $t_{m-1}$ such that $e <_{\pi(t_{m-1})} e_k$. Every such edge
%also satisfies $e<_{\pi(T)} e_k$ because the order $\pi(t_{m-1})$ is
%consistent with $\pi(T)$. The lemma follows.
%\end{proof}

\begin{lemma}\label{order}
In the above construction, $p_0 <_{\pi(T)} \dots <_{\pi(T)} p_n$.
\end{lemma}

\begin{proof}
Since $\Pi(G)$ is a proper set of tree orders, it follows that the root $0$ is
the smallest vertex of $T$ in the order $\pi(T)$. Hence $p_0 <_{\pi(T)} p_1$.
Suppose that $p_0 <_{\pi(T)} \dots <_{\pi(T)} p_{m}$ for some $1\leq
m \leq n-1$. We show that $p_{m} <_{\pi(T)} p_{m+1}$. 
We consider the
following 
two cases.

\vskip5pt
Case 1: $p_{m+1} \notin V_{m}$. Then the number of edges from
$p_{m+1}$ to $t_{m-1}$ is at most $b_{p_{m+1}}$. Since $p_{m+1} \in
V_{m+1}$, the number of edges from $p_{m+1}$ to $t_m$ is at least
$b_{p_{m+1}}+1$. It follows that there is at least one edge
$(p_{m+1},p_m)$ in $G$ and that $p_{m+1}$ is adjoined to $t_m$ by
means of such an edge. Thus $p_m <_{\pi(T)} p_{m+1}$ because $\Pi(G)$ is
a proper set of tree orders.

\vskip5pt
Case 2: $p_{m+1} \in V_{m}$. Let $e_{p_{m+1}}$ be the edge from
$p_{m+1}$ to $t_{m-1}$ such that exactly $b_{p_{m+1}}$ edges $e$ from $p_{m+1}$ to
$t_{m-1}$ satisfy $e <_{\pi(t_{m-1})} e_{p_{m+1}}$. 
Since $p_m$ is the
largest vertex of $t_m$ in the order $\pi(T)$ and hence in the order
$\pi(t_m)$, and $e_{p_{m+1}}$ goes from $p_{m+1}$
to $t_{m-1}=t_m-p_m$,
it follows that $e <_{\pi(t_{m-1})} e_{p_{m+1}}$ if and
only if $e<_{\pi(t_m)} e_{p_{m+1}}$ because the order $\pi(t_{m-1})$
is consistent with $\pi(t_m)$. 
Therefore, exactly $b_j$ edges $e$ from $p_{m+1}$ to $t_m$
satisfy $e <_{\pi(t_m)} e_{p_{m+1}}$, hence $p_{m+1}$ is adjoined to
$t_m$ by means of the edge $e_{p_{m+1}}$.

Let $e_{p_m}$ be the edge of $T$ coming out of $p_m$, and let $t$ be
the tree in the construction of $T$ obtained by adjoining the vertices
of $V_m$ to $t_{m-1}$. Let $t'$ be the tree obtained from $t_{m-1}$ by
adjoining the vertices $p_m$ and $p_{m+1}$ by means of the edges
$e_{p_m}$ and $e_{p_{m+1}}$. Then $t'$ is a subtree of both $t$ and
$T$. By choice of $p_m$, we have $p_m <_{\pi(t)} p_{m+1}$, so $p_m
<_{\pi(t')} p_{m+1}$ and $p_m <_{\pi(T)} p_{m+1}$ because the order
$\pi(t')$ is consistent with both $\pi(t)$ and $\pi(T)$.
\end{proof}

We now check that $\Theta$ and $\Phi$ are inverses of each other.

\begin{lemma}\label{xiphi}
$\Theta(\Phi(P)) = P$ for $P\in\mathcal P_G$.
\end{lemma}

\begin{proof}
Put $P = (b_1,\dots,b_n)$ and $T = \Phi(P)$.
Consider the process of constructing $T$.
For $j\in\{1,\dots,n\}$, we have $j=p_m$ for some $1\leq m \leq n$. 
Let $e_j$ be the edge of $T$ coming out of $j$. The edge $e_j$ goes
from $j$ to $t_{m-1}$. Since the set of vertices of $t_{m-1}$ is
$\{p_0,\dots,p_{m-1}\}$, it follows from Lemma \ref{order} that if an
edge
$e$ coming out of $j$ satisfies $e <_{\pi(T)} e_j$, then $e$ goes from
$j$ to $t_{m-1}$. Thus, $e <_{\pi(T)} e_j$ if
and only if $e <_{\pi(t_{m-1})} e_j$ because the order $\pi(t_{m-1})$
is consistent with $\pi(T)$. By construction of $T$,
the number of edges $e$ satisfying $e <_{\pi(t_{m-1})} e_j$ is $b_j$,
hence the number of edges $e$ satisfying $e <_{\pi(T)} e_j$ is
also $b_j$. We conclude that $\Theta(T)=P$.
\end{proof}

\begin{lemma}\label{phixi}
$\Phi(\Theta(T')) = T'$ for $T'\in\mathcal T_G$.
\end{lemma}

\begin{proof}
Put $P=\Theta(T')=(b_1,\dots,b_n)$. Consider the process of constructing
$T = \Phi(P)$. We show by induction that for $0\leq m \leq n$, the
tree $t_m$ is a subtree of $T'$ and that $p_0,\dots,p_m$ are the
smallest $m+1$ vertices in the order $\pi(T')$. 
Since the root $0$ is the smallest vertex in
$\pi(T')$, the assertion is true for $m=0$.

Now, suppose that $t_{m-1}$ is a subtree of $T'$ and that $p_0, \dots,
p_{m-1}$ are the smallest $m$ vertices in the order $\pi(T')$. Let $k$
be the $(m+1)$-th smallest vertex in the order $\pi(T')$, and let
$e'_k=(k,i)$ be the edge coming out of $k$ in $T'$. Then $i
<_{\pi(T')} k$, so $i\in
\{p_0,\dots,p_{m-1}\}$ and $i\in t_{m-1}$. Hence if an edge $e$ coming
out of $k$ satisfies $e <_{\pi(T')} e'_k$, then $e$ goes from $k$ to
$t_{m-1}$. There are $b_k$ edges $e$ satisfying $e <_{\pi(T')}
e'_k$. These $b_k$ edges together with the edge $e'_k$ give $b_k+1$
edges going from $k$ to $t_{m-1}$. It follows that $k\in V_m$.

As before, for every $j\in V_m$, let $e_j$ be the edge from $j$ to
$t_{m-1}$ such that exactly $b_j$ edges $e$ from $j$ to $t_{m-1}$
satisfy $e <_{\pi(t_{m-1})} e_j$. 
Since the vertices of $t_{m-1}$ are the smallest $m$ vertices in the
order $\pi(T')$, it follows that if an edge $e$ coming out of $j$
satisfies $e <_{\pi(T')} e_j$, then $e$ goes from $j$ to
$t_{m-1}$. Thus, 
$e <_{\pi(T')} e_j$ if
and only if $e <_{\pi(t_{m-1})} e_j$ because the order $\pi(t_{m-1})$
is consistent with $\pi(T')$. There are $b_j$ edges $e$ satisfying $e
<_{\pi(t_{m-1})} e_j$, hence there are $b_j$ edges $e$ satisfying $e
<_{\pi(T')} e_j$. It follows from the choice of $b_j$ that $e_j$ is an
edge of $T'$. Therefore, the tree $t$ obtained by adjoining the
vertices $j\in V_m$ by means of the edges $e_j$ is a subtree of
$T'$. Consequently, the smallest vertex $p_m$ of $V_m$ in the order $\pi(t)$
is the smallest vertex of $V_m$ in the order $\pi(T')$. Since $k$
is the smallest vertex of $U_m$ in the order $\pi(T')$ and
$k\in V_m\subseteq U_m$, it follows that $p_m=k$. The induction step
is complete.

Finally, we obtain $T' = t_n = T$.
\end{proof}

Theorem \ref{main} follows from Lemmas \ref{xiphi} and \ref{phixi}.
\end{proof}

\section{Examples}\label{examples}

In this section we give examples of proper sets of tree orders and the
resulting bijections between $\mathcal T_G$ and $\mathcal P_G$
from the family of bijections defined in Section \ref{bijections}.

We begin by introducing the \emph{breadth-first search order} $\pi_{bf}(T)$
on the
vertices of a tree $T\in\mathbb T_G$. For a vertex $i\in T$, we
define the \emph{height} $h_T(i)$ of $i$ in $T$ to be the number of
edges in the unique path from $i$ to the root $0$. We set $i
<_{\pi_{bf}(T)} j$, or $i <_{bf} j$, if $h_T(i) < h_T(j)$ or else if
$h_T(i)=h_T(j)$ and $i<j$. It is easy to check that $\pi_{bf}(T)$ is a
total order on the vertices of $T$ and that $\Pi_{bf}(G) =
\{\pi_{bf}(T)\ |\ T\in \mathbb T_G\}$ is a proper set of tree orders.

The \emph{depth-first search order} $\pi_{df}(T)$ on the vertices of
  a tree $T\in\mathbb T_G$ is defined as follows. For a vertex $i\in
  T$, let $T(i)$ denote the branch of $T$ rooted at $i$. In other
  words, $T(i)$ consists of all vertices $k$ of $T$ such that the
  unique path
  from $k$ to $0$ in $T$ contains $i$. If $(i,\ell)$ is an edge of
  $T$, then we set $\ell <_{\pi_{df}(T)} i$.
Furthermore, if $(j,\ell)$
  is an edge of $T$ such that $i<j$, 
then we set $i' <_{\pi_{df}(T)} j'$ for $i'\in
  T(i)$ and $j'\in T(j)$. We use the symbol $<_{df}$ with the same
  meaning as $<_{\pi_{df}(T)}$. It is not hard to see that
  $\Pi_{df}(G)=\{\pi_{df}(T)\ |\ T\in\mathbb T_G\}$ is a proper set of
  tree orders.

Our third example is the \emph{vertex-adding order} $\pi_{va}(T)$ on
the vertices of $T\in\mathbb T_G$. Construct the sequence
$p_0,\dots,p_{|T|-1}$ inductively as follows. Set $p_0=0$, and, for $1\leq m
\leq |T|-1$, let $p_m$ be
the smallest vertex $j$ in $G-\{p_0,\dots,p_{m-1}\}$ such that there is an
edge in $G$ from $j$ to $\{p_0,\dots,p_{m-1}\}$. Note that
the sequence $p_0,\dots,p_{|T|-1}$ contains each vertex of $T$ exactly once.
Put $p_0
<_{\pi_{va}(T)} \dots <_{\pi_{va}(T)} p_{|T|-1}$, and let the symbol $<_{va}$
have the same meaning as $<_{\pi_{va}(T)}$. Clearly, $<_{va}$ is a
total order on the vertices of $T$. Also, $\Pi_{va}(G) =
\{\pi_{va}(T)\ |\ T\in\mathbb T_G\}$ is a proper set of tree
orders. Indeed, if $t$ is a subtree of $T$, then adding or not adding
a vertex of
$T-t$ to $\{p_0,\dots,p_{m-1}\}$ does not affect the order in which
the vertices of $t$ are added.

\begin{figure}
\begin{center}
\input{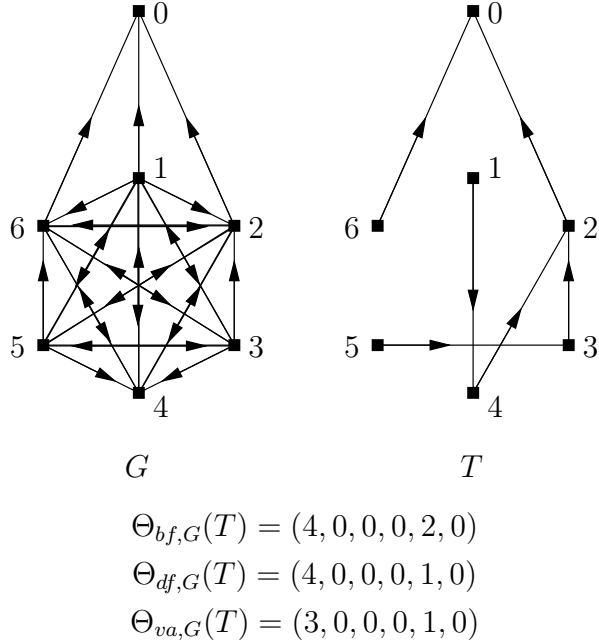}
\end{center}
\caption{An example of constructing $\Theta_{bf,G}(T)$,
  $\Theta_{df,G}(T)$, and $\Theta_{va,G}(T)$.}\label{exempel}
\end{figure}

Let $\Theta_{bf,G}$, $\Theta_{df,G}$, and $\Theta_{va,G}$ be the maps
$\Theta_{\Pi_{bf},G}$, $\Theta_{\Pi_{df},G}$, and
$\Theta_{\Pi_{va},G}$ constructed in Section \ref{bijections}. Figure
\ref{exempel} shows a sample graph $G$ and a spanning tree
$T\in\mathcal T_G$. To compute $\Theta_{bf,G}(T)$,
  $\Theta_{df,G}(T)$, and $\Theta_{va,G}(T)$, we first determine the
orders $\pi_{bf}(T)$, $\pi_{df}(T)$, and $\pi_{va}(T)$. 
We have $h_T(0) = 0$, $h_T(2) = h_T(6) = 1$, $h_T(3) = h_T(4) =2$,
and $h_T(1) = h_T(5) = 3$, so
$$
0 <_{bf} 2 <_{bf} 6 <_{bf} 3 <_{bf} 4 <_{bf} 1 <_{bf} 5.
$$

Next, we determine $\pi_{df}(T)$. Applying the depth-first search rule
with $\ell=0$, we get $0 <_{df} \{1,2,3,4,5\} <_{df} 6$ because $T(2)$
contains vertices $1$, $2$, $3$, $4$, and $5$, and $T(6)$ contains a
single vertex $6$. Taking $\ell = 2$ we get $0 <_{df} 2 <_{df} \{3,5\}
<_{df} \{1,4\} <_{df} 6$. Finally, taking $\ell = 3$ and $\ell = 4$ we
get
$$
0 <_{df} 2 <_{df} 3 <_{df} 5 <_{df} 4 <_{df} 1 <_{df} 6.
$$

Also, the order $\pi_{va}(T)$ is the following:
$$
0 <_{va} 2 <_{va} 3 <_{va} 4 <_{va} 1 <_{va} 5 <_{va} 6.
$$

The edge coming out of vertex $1$ in $T$ is $(1,4)$. The relation $e
<_{bf} (1,4)$ is satisfied for $e = (1,0),
(1,2), (1,3), (1,6)$, so the first component of $\Theta_{bf,G}(T)$ is
$4$. Similarly, $e <_{df} (1,4)$ holds for $e = (1,0), (1,2),
(1,3), (1,5)$, so the first component of $\Theta_{df,G}(T)$ is $4$. 
The relation $e <_{va} (1,4)$ holds for $e = (1,0), (1,2), (1,3)$, so
the first component of $\Theta_{va,G}(T)$ is $3$. 
The other components are computed in the same way.
Figure \ref{exempel}
shows the values of $\Theta_{bf,G}(T)$, $\Theta_{df,G}(T)$, and
$\Theta_{va,G}(T)$.

\begin{figure}                                                                 \
                                                                                
\begin{center}                                                                 \
                                                                                
\input{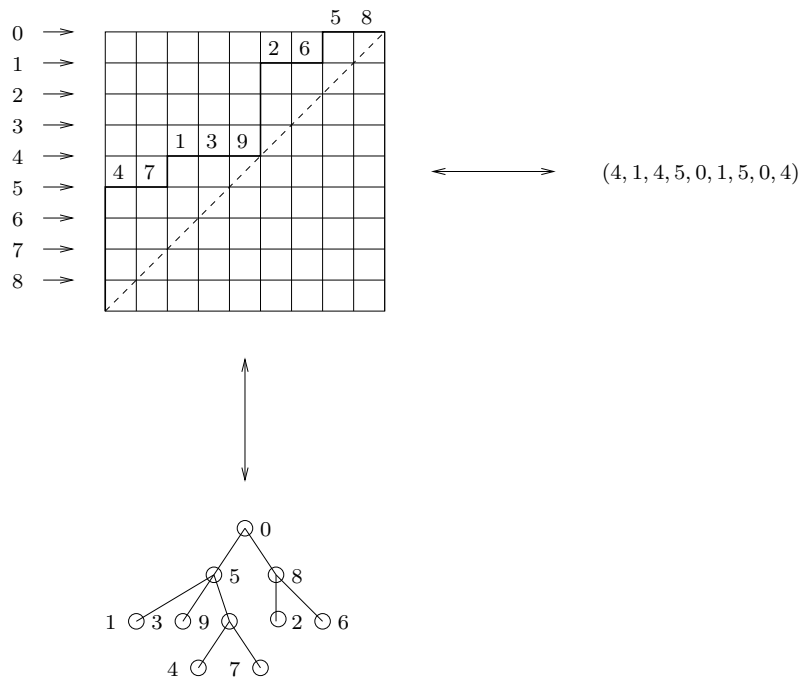}                                                     \
                                                                                
\end{center}                                                                   \
                                                                                
\caption{A bijection between trees and parking functions involving              
labeled Dyck paths.}\label{postnikovs}                                          
\end{figure}

Note that for $G=K_{n+1}$, the presented construction yields a                 
family of bijections between the classical parking functions and trees         
on $n+1$ labeled vertices. This family includes some of the well-known         
bijections. For example, using the vertex-adding tree order results in         
the following simple correspondence defined in terms of drivers and            
parking spots: given a parking function $(b_1,\dots,b_n)$, the                 
corresponding tree is obtained by introducing the edge $(i,j)$                 
whenever driver $j$ ended up parking in spot $b_i-1$, and the edge             
$(i,0)$ whenever $b_i=0$.

Another bijection 
involving labeled Dyck paths as an intermediate object,
communicated to us by A. Postnikov, 
results if the right-to-left depth first search tree order is
used (this order is the same as the depth first search order described
above except that larger numbers are given priority among the children
of the same vertex). Given a parking function $P$, we write numbers 
$1$ through $n$ in the $n\times n$ square so that all numbers $j$ such
that $b_j=i$ appear in the $i$-th row in increasing order, and the
numbers
in a lower row appear to the left of the numbers in a higher row. Such
an arrangement defines a Dyck path from the lower-left corner to the
upper-right corner of the square, with horizontal steps labeled with
integers between $1$ and $n$; see Figure \ref{postnikovs}. To get the
spanning
tree $T$ corresponding to $P$, start from the upper-right corner of
the square and proceed to the lower-left corner along the Dyck path,
keeping track of the current vertex, initially set to be $0$. At each
horizontal step labeled $i$, connect the vertex $i$ to the current
vertex, and at each vertical step, replace the current vertex with its
successor in the right-to-left depth first seach order on the tree
constructed so far. It is not hard to show that the obtained tree $T$
is precisely $\phi(P)$ for the right-to-left depth first search tree order.

The bijection obtained using the breadth first search tree order is discussed in Section \ref{sandpiles} in connection with the sandpile model.

\section{More proper sets of tree orders}\label{more}

We now present a method for constructing proper
sets of tree orders. Let $\langle \sigma_1,\dots,\sigma_\ell \rangle$
denote the path consisting of the edges
$(\sigma_\ell,\sigma_{\ell-1})$, $(\sigma_{\ell-1},\sigma_{\ell-2})$,
\dots, $(\sigma_1,0)$. Also, let $\langle \rangle$ denote the path
consisting of the vertex $0$ alone.
Define $\mathbb A_G$ to be the set of paths $\langle
\sigma_1,\dots,\sigma_\ell\rangle$ in $G$ such that
$\sigma_1,\dots,\sigma_\ell$ are distinct vertices of $G-\{0\}$, where
$\ell \geq 0$.
Let $\prec$ be a partial order on $\mathbb A_G$ satisfying the following
conditions:
\begin{enumerate}

\item[(i)] if $A\cap A' \in \mathbb A_G$ for some $A,A'\in\mathbb
  A_G$, then $A$ and $A'$ are comparable;

\item[(ii)] $\langle \sigma_1,\dots,\sigma_{\ell'}\rangle \prec
  \langle \sigma_1,\dots,\sigma_{\ell'},\dots,\sigma_\ell\rangle$ for
  $\ell'<\ell$.

\end{enumerate}

For a tree $T\in\mathbb T_G$ and a vertex $i\in
T$, let $A_T(i)\in \mathbb A_G$ 
be the unique path in $T$ from $i$ to $0$. Introduce the order
$\pi_\prec(T)$ on the vertices of $T$ in which $i <_{\pi_\prec(T)} j$
whenever $A_T(i) \prec A_T(j)$. Put $\Pi_\prec(G) = \{\pi_\prec(T)\ |\
T\in\mathbb T_G\}$.

\begin{proposition}
$\Pi_\prec(G)$ is a proper set of tree orders.
\end{proposition}

\begin{proof}
Let $T\in\mathbb T_G$, and let $i$ and $j$ be vertices of
$T-\{0\}$. Since $A_T(i)$ and $A_T(j)$ are the unique paths in $T$
from $i$ and
$j$ to $0$, it follows that $A_T(i)\cap A_T(j)\in \mathbb
A_G$. Therefore, $\pi_\prec(T)$ is a total order on the vertices of
$T$, by property (i) of $\prec$.

If $(j,i)$ is an edge of $T$, then $A_T(i) = \langle
\sigma_1,\dots,\sigma_\ell,i\rangle$ and $A_T(j) = \langle
\sigma_1,\dots,\sigma_\ell, i, j\rangle$, so $A_T(i) \prec A_T(j)$, by
property (ii) of $\prec$, so $i <_{\pi_\prec(T)} j$. 

If $t$ is a subtree of $T$, then $A_t(i) = A_T(i)$ for all vertices $i\in t$, 
so the order $\pi_\prec(t)$ is consistent with the order
$\pi_\prec(T)$.

The proposition follows.
\end{proof}

\begin{figure}                                                                  
\begin{center}                                                                 
\scalebox{0.6}{ 
\input{exempel3.pstex_t}                                                        
}
\end{center}                                                                    
\caption{}\label{exempel3}                                                      
\end{figure}       

The orders $\pi_{bf}(T)$, $\pi_{df}(T)$, and $\pi_{va}(T)$ described
in Section \ref{examples} can be
obtained as $\pi_\prec(T)$ via an appropriate choice of
$\prec$. Setting $\prec$ to be the lexicographic order on the paths
$\langle \sigma_1,\dots,\sigma_\ell\rangle$ viewed as sequences of
integers yields the order $\pi_{df}(T)$. To obtain $\pi_{bf}(T)$, set
$\langle \sigma_1,\dots,\sigma_\ell \rangle \prec \langle
\sigma'_1,\dots,\sigma'_{\ell'}\rangle$ if $\ell<\ell'$, or else if
$\ell=\ell'$ and $\sigma_\ell<\sigma'_{\ell'}$. Finally, setting
$\prec$ to be the order in which $A \prec A'$ whenever 
$A\cap A'\in\mathbb A_G$, and
the largest
vertex of $A\backslash A'$ is smaller than the largest vertex of
$A'\backslash A$, yields the order $\pi_{va}(T)$.

We can obtain other proper sets of tree orders from partial
orders $\prec$ on $\mathbb A_G$ satisfying the conditions above. For
example, we can set $A=\langle \sigma_1,\dots,\sigma_\ell\rangle \prec
A'=\langle \sigma'_1,\dots,\sigma'_{\ell'}\rangle$ whenever the
increasing rearrangement of $A$ is smaller than that of $A'$ in the
lexicographic order. Another example is setting $A\prec A'$ if
$\sum \sigma_k < \sum \sigma'_k$, or else if $\sum \sigma_k =
\sum \sigma'_k$ and $\sigma_\ell<\sigma'_{\ell'}$.

Similar examples of partial orders 
on $\mathbb{A}_G$
yielding proper sets of tree orders
can be obtained by using an arbitrary numbering of the edges of $G$ instead 
of vertex labels.

%\begin{figure}
%\begin{center}
%\input{exempel3.pstex_t}
%\end{center}
%\caption{}\label{exempel3}
%\end{figure}

It is worth noting that not all proper sets of tree orders are induced by
a partial order on $\mathbb{A}_G$ satisfying the above
conditions. Consider the following simple example. Let $G$ be the
graph shown in Figure \ref{exempel3}. Let $e_{1,2}$ be the two edges
of $G$ going from vertex $3$ to vertex $1$, and let $f_{1,2}$ be the
two edges going from vertex $4$ to vertex $2$. For $1\leq i,j\leq 2$,
let $T_{ij}$ be the spanning tree of $G$ containing edges $e_i$ and
$f_j$. Let $\Pi=\{\pi(T_{ij})\ |\ 1\leq i,j\leq 2\}$ 
be the proper set of tree orders defined as follows:
$$
0 <_{\pi(T_{ij})} 1   <_{\pi(T_{ij})} 2 <_{\pi(T_{ij})} 3
<_{\pi(T_{ij})} 4
$$
for $i\neq j$, and
$$
0 <_{\pi(T_{ii})} 1   <_{\pi(T_{ii})} 2 <_{\pi(T_{ii})} 4                       
<_{\pi(T_{ii})} 3.  
$$
Let $A_{e_i}$ (resp.\ $A_{f_i}$) be the unique path in
$\mathbb{A}_G$
from
vertex $3$ (resp.\ $4$) to the root $0$ containing the edge $e_i$
(resp.\ $f_i$). Then in order for $\Pi$ to be induced by some partial
order $\prec$ on $\mathbb{A}_G$, we must have $A_{e_1} \prec A_{f_2}$
so that relation $3 <_{\pi(T_{12})} 4$ holds. Similarly, to achieve
relations $4 <_{\pi(T_{22})} 3$, $3 <_{\pi(T_{21})} 4$, and $4
<_{\pi(T_{11})} 3$, we must have $A_{f_2} \prec A_{e_2}$, $A_{e_2}
\prec A_{f_1}$, and $A_{f_1}\prec A_{e_1}$. We obtain a contradiction
$A_{e_1} \prec A_{e_1}$, hence $\Pi$ is not induced by a partial order 
on $\mathbb{A}_G$.

\section{$G$-parking functions and the sandpile model}\label{sandpiles}

In \cite{Cori}, Cori and Le Borgne construct a family of bijections
between the rooted spanning trees of a digraph $G$ and the recurrent
states of the sandpile model defined on $G$. It was shown by
Gabrielov \cite{Gab1} that if
for all vertices of $G$ except the root, the out-degree is greater than or equal to the
in-degree, then recurrent states coincide with the so called \emph{allowed
configurations} of the model, which 
correspond to $G$-parking
functions: if $d_i$ is the out-degree of vertex $i$, then $(u_1,\dots,u_n)$ is an allowed configuration
if and only if $(d_1-u_1,\dots,d_n-u_n)$ is a
$G$-parking function.
In particular, this observation is valid for symmetric
graphs, in which the number of edges from $i$ to $j$ is equal to the
number of edges from $j$ to $i$ for all $i\neq j$; such graphs can be
naturally 
viewed as undirected graphs. Thus for these graphs the result of
Cori and Le Borgne provides a bijective correspondence between
rooted spanning trees of $G$ and $G$-parking functions. For the rest
of the section, we assume that $G$ is a symmetric graph.

The construction described in \cite{Cori} begins by fixing an
arbitrary order on the edges of $G$. 
Given a spanning tree $T$ of $G$, an edge $e$ in $G-T$ is called
\emph{externally active} 
with respect to $T$ if in the unique cycle of $T+e$,
the edge $e$ is the smallest in the chosen order. 
A key property of the obtained
bijection is that the sum of the values of a recurrent state 
is equal to the number of 
externally active edges
with respect to the corresponding spanning tree. It follows that in
the resulting bijection between $G$-parking functions and spanning
trees, $G$-parking functions with the same sum of values are mapped to
spanning trees with the same number of externally active edges.

To show that the bijections presented in this paper are substantially
different from the ones in \cite{Cori}, consider the case $G=K_{n+1}$,
and let $P$ be the path obtained as follows: start at the root vertex
$0$, and then append the remaining vertices one by one, so that at
each step the appended edge is the smallest, in the chosen edge order,
among all edges that can possibly be appended. There are no externally
active edges with respect to $P$ since every edge $(i,j)$ not in
$P$, where $j$ is closer to the root in $P$ than $i$, is greater than
the edge $(i',j)$, where $i'$ is the vertex appended after $j$ in the
construction of $P$, by choice of $i'$. On the other hand, if a path $P'$
does not include the smallest edge in the chosen edge order, then this
edge is externally active with respect to $P'$. Hence there is a
different number of externally active edges with respect to $P$ and
$P'$. However, every bijection $\Theta_{\Pi,G}$ maps both $P$ and $P'$
to permutations of $(0,\dots,n-1)$, so the sum of values of the
corresponding $G$-parking functions is the same. Hence for
$G=K_{n+1}$, none of the bijections $\Theta_{\Pi,G}$ coincides with
a bijection from the family constructed in \cite{Cori}.

Dhar defined the {\it {burning algorithm}} for determining                                                                  
whether a given configuration is allowed; see \cite{Znam}. 
In our setting this task                                                              
corresponds to the question whether a function $P\ :\ \{1,\dots,n\} \rightarrow \mathbb{N}$ is a $G$-parking                
function, and an equivalent formulation of Dhar's burning algorithm is the                                                  
following. We mark vertices of the graph, starting with the root $0$.                                                       
At each iteration of the algorithm, we                                                                                      
mark all vertices $v$ that have more marked neighbors                                                                       
than the value of the function at $v$. If in the end                                                                        
all vertices are marked, then we have a $G$-parking function, as it                                                         
is not hard to see directly from definition. Conversely, for every                                                          
$G$-parking function, this algorithm marks all vertices.                                                                    
                                                                                                                            
We claim that our bijection corresponding to the breadth first search                                                       
order $\pi_{bf}$ is a natural generalization of Dhar's                                                                      
algorithm. Given a parking function $P=(b_1,\dots,b_n)$,                                                                    
perform the construction of $T=\Phi_{\Pi, G}(P)$ as                                                                         
described above. We know that $T$                                                                                           
contains all verices if and only if we started with a $G$-parking                                                           
function. Let us group the vertices of $T$ by height, setting $W_i$ to                                                      
be                                                                                                                          
the set of vertices of $T$ of height $i$.                                                                                   
                                                                                                                            
\begin{proposition}\label{burningbreadthfirst}                                                                              
$W_i$          
is exactly the set of vertices marked at the $i$-th step of the burning algorithm.
\end{proposition}                                                                                                           
                                                                                                                            
\begin{proof}
For $i=0$ the claim is true because the root $0$ is marked                                                                  
at the $0$-th step. We prove the claim by induction.                                                                        
Suppose that for $k<i$, the vertices in $W_k$ are marked at the                                                             
$k$-th step of the Dhar's algorithm.                                                                                        
Let $e_j=(j,w_j)$ be the edge going out of $j$ in $T$.                                                                      
Each vertex $j\in W_i$ has more than $b_j$ edges                                                                                  
going to vertices not larger than $w_j$ in                                                                                  
$\pi_{bf}$ order. All vertices not                                                                                          
larger than $w_j$ are in $\cup_{k<i} W_k$ since $w_j\in W_{i-1}$.                                                           
%Thus                                                                                                                       
%vertices from $V_i$ have no influence on ability to add to the tree a                                                      
%vertex from $V_i$.                                                                                                         
Therefore, all vertices in $W_i$ are marked                                                                                  
at the $i$-th step of Dhar's algorithm. On the other hand, every vertex                                                     
marked at the $i$-th step of the algorithm in our                                                                          
is to be attached in $T$ to a vertex from $\cup_{k<i} W_k$ since we                                                       
add vertices to $T$ in the order $\pi_{bf}$. Thus                                                                           
each such vertex is in $W_i$. Hence Dhar's burning                                                                          
algorithm is realized by our bijection for the breadth first                                                                
search tree order.                                            
\end{proof}

\section{Acknowledgments}\label{thanks}

The authors would like to thank Prof. Alexander Postnikov for suggesting this
problem and for many helpful discussions. We are also grateful to Prof. Igor
Pak and Prof. Richard Stanley for feedback and advice.


\begin {thebibliography}{[L]}

\bibitem {Cori}
R. Cori and Y. Le Borgne: The sand-pile model and Tutte polynomials,
{\it Advances in Applied Mathematics} {\bf 30} (2003), 44--52.

\bibitem {Dhar}
D. Dhar: Self-organised critical state of the sandpile automaton models, 
{\it Physical Review Letters} {\bf 64} (1990), no. 14, 1613--1616.

\bibitem {Fran}
J. Francon: Acyclic and parking functions, {\it J. Combinatorial
  Theory, Ser. A} {\bf 18} (1975), 27--35.

\bibitem {Gab1}
A. Gabrielov: Abelian avalanches and Tutte polynomials, {\it Physica A}
{\bf 195} (1993), 253--274.

\bibitem {Gab}
A. Gabrielov: Asymmetric abelian avalanches and sandpiles, preprint 93--65, 
MSI, Cornell University, 1993.

\bibitem {Iva}
E. V. Ivashkevich, V. B. Priezzhev: Introduction to the sandpile model, {\it 
Physica} A {\bf 254} (1998), 97--116.

\bibitem{Znam}                                                                  
R. Meester, F. Redig, and D. Znamenski:                                         
The Abelian sandpile; a mathematical introduction, {\it Markov                  
  Processes and Related Fields} {\bf 7} (2001), 509--523.

\bibitem {Post}
A. Postnikov, B. Shapiro: Trees, parking functions, syzygies, and 
deformations of monomial ideals, preprint arXiv:math.CO/0301110, 2003.

\bibitem {Stan}
R. P. Stanley: {\it Enumerative Combinatorics, Volume 2}, Cambridge Studies 
in Advanced Mathematics {\bf 62}, Cambridge University Press, Cambridge, 
1999.

%\bibitem{Znam}
%R. Meester, F. Redig, and D. Znamenski: 
%The Abelian sandpile; a mathematical introduction, {\it Markov
%  Processes and Related Fields} {\bf 7} (2001), 509--523.

\end {thebibliography}
\end {document}